\documentclass[conference,10pt,twoside]{IEEEtranTMBMC}
\normalsize
\usepackage[cmex10]{amsmath}
\usepackage{amsfonts}
\usepackage{amssymb}
\usepackage{array}
\usepackage{cite}
\usepackage{url}
\usepackage{color}

\usepackage{fixltx2e}
\usepackage[pdftex]{graphicx}

\title{Inferring biological networks by sparse identification of nonlinear dynamics}
\author{
	Niall~M. Mangan\IEEEauthorrefmark{1}\IEEEauthorrefmark{3}, 
	Steven L.~Brunton\IEEEauthorrefmark{2}, 
	Joshua L.~Proctor\IEEEauthorrefmark{3}, 
	and~J. Nathan~Kutz\IEEEauthorrefmark{1}
	\\
	 \IEEEauthorrefmark{1}{Department of Applied Mathematics, University of Washington, Seattle, WA 98195 USA }
	\\ \IEEEauthorrefmark{2}{Department of Mechanical Engineering and the eScience Institute, University of Washington, Seattle, WA 98195 USA}
	 \\ \IEEEauthorrefmark{3}{Institute for Disease Modeling, 3150 139th Ave SE
Bellevue, WA 98005}
	}
	
\begin{document}

\maketitle

\begin{abstract}
Inferring the structure and dynamics of network models is critical to understanding the functionality and control of complex systems, such as metabolic and regulatory biological networks. 
The increasing quality and quantity of experimental data enable statistical approaches based on information theory for model selection and goodness-of-fit metrics. 
We propose an alternative method to infer networked nonlinear dynamical systems by using sparsity-promoting $\ell_1$ optimization to select a subset of nonlinear interactions representing dynamics on a fully connected network. 
Our method generalizes the sparse identification of nonlinear dynamics (SINDy) algorithm to dynamical systems with rational function nonlinearities, such as biological networks.
%
%
We show that dynamical systems with rational nonlinearities may be cast in an implicit form, where the equations may be identified in the null-space of a library of mixed nonlinearities including the state and derivative terms; this approach applies more generally to implicit dynamical systems beyond those containing rational nonlinearities. 
This method, implicit-SINDy, succeeds in inferring three canonical biological models: Michaelis-Menten enzyme kinetics, the regulatory network for competence in bacteria, and the metabolic network for yeast glycolysis.
\end{abstract}


\IEEEpeerreviewmaketitle

\section{Introduction}
\IEEEPARstart{N}{etwork} science is of growing and critical importance across the physical, biological and engineering sciences.  In biology, both the quality and quantity of modern data has inspired new mathematical techniques for inferring the complex interactions and connections between nodes in metabolic and regulatory networks. Discovering the connectivity and structure of such networks is critical in understanding the functionality and control decisions enacted by the network in tasks such as cell differentiation, cell death, or directing metabolic flux.  Methods based on information theory provide rigorous statistical criteria for such model selection and network inference tasks.   
For example, partnering the Kullback-Leibler (KL) divergence \cite{kullback1951,kullback1959}, a measure of information loss between empirically collected data and model generated data, and the Akaike information criteria (AIC)  \cite{akaike1973,akaike1974}, a relative estimate of information loss across models balancing model complexity and goodness-of-fit, allows for a principled model selection criteria \cite{burnham2002}.  
However, for nonlinear dynamical networks, such information theoretic approaches face severe challenges since the number of connections and possible functional forms lead to a combinatorially large number of models to be evaluated with the KL/AIC framework.  

We propose an alternative method to infer the dynamics and connectivity of biological networks motivated by machine learning methods including overcomplete libraries~\cite{bishop2006pattern,murphy2012machine,James2013book} and sparse  regression~\cite{Tibshirani:1996, Zou2005}.  
We generalize the {\em sparse identification of nonlinear dynamics} (SINDy)~\cite{Brunton2016pnas} algorithm to an implicit formulation that allows the dynamics to contain rational functions.
We demonstrate the accuracy and robustness of the method, {\em implicit-SINDy}, on three representative biological models, showing that our approach gives a compelling alternative to information-theoretic methods of model selection.

\subsection{Biological networks}

Biological networks produce a diverse range of functional activities.  Regulatory and metabolic networks are critical for cellular function.  Breakdown of the function and control circuits of such networks can lead to cancer and other deadly diseases, motivating attempts to control of such networks using the tools of genetic engineering and pharmacology. 
%
%
The network dynamics in these systems can often be modeled using mass-action kinetics, producing a relatively constrained set of  network motifs \cite{Alon2007a, Alon2007}.  Generally, biological regulatory and metabolic networks are considered sparse and sparsity has been used as a criteria for inferring linear network models \cite{Grigorov2005, Andrecut2008, Oksuz2013, Nordling2011}. 
However, accurately and robustly characterizing such biological networks, both in terms of their unique connectivity structure and functional dynamics, remains an extremely challenging task due to the underlying nonlinear dynamics.

With the emergence of  large amounts of high-quality experimental data\cite{Spiller2010, Wu2012}, new  opportunities exist for data-driven mathematical modeling of these biological networks.  Rapid, robust and accurate model identification can greatly accelerate our understanding and control of critical biological network functions so that disease treatment and metabolic engineering protocols can be proposed.  Indeed, the rich, dynamic data on such biological networks~\cite{Spiller2010, Wu2012} allows for greatly improved statistical and machine learning methods that can help characterize the fundamental dynamic interactions in complex, biological networks.

\subsection{Model selection and information theory}
Biological network models can be developed using principled dynamical systems techniques based upon conjectured relationships between variables (nodes).  Such models are constructed and iteratively modified to be consistent with experiment.  Historically, such trial-and-error development could take many years before a model is quantitatively predictive.   Statistical methods advocate techniques where information criteria for model selection assigns a numerical value for each candidate model.  The models with the lowest values correspond to the least information loss with respect to all candidate models considered.  Common examples of information criteria are the AIC~ \cite{akaike1973} and the Bayesian information criteria (BIC)~\cite{Bayesian1978}.  See \cite{Hjorth2008} for a review of modern usage of information criteria.  

Unfortunately, information criteria rely on the practitioner positing a set of reasonable models. If the network interactions and functional nonlinearities are unknown, then it is unlikely that one can guess the correct network structure and nonlinearities in the system.  Of course, one may select from a large space of trial models, but each individual model would need to be simulated, compared to the data, and given a numerical information criteria score.  This can be computationally intractable as the number of nonlinear interactions between variables scales combinatorially with the number of variables. 

To highlight the complexity of model selection, consider selecting a specific model from all possible polynomials of degree $d$ in $n$ variables.  For $n$ variables, the number of possible monomials with degree less than or equal to $d$ is given by
\begin{eqnarray}
N_m = {n + d \choose d}.
\end{eqnarray}
The number of possible polynomial structures that may be formed by assigning nonzero coefficients to these $N_m$ monomials is given by:
\begin{eqnarray}
N_p = \sum_{k=1}^{N_m} { N_m \choose k}.\label{eq:nummodels}
\end{eqnarray}
The number of possible polynomial structures $N_p$ may be thought of as summing over all possible polynomials with only $k$  monomials.   
For example, consider polynomials up to degree $d=4$ with $n=5$ variables, as is required for the metabolic network example in this paper.  
This leads to $N_m=126$ and $N_p\approx 10^{38}$.  
The number of possible rational functions given general numerator and denominator polynomials is even larger. 
Simulating, evaluating, and comparing this many models in an information criteria framework is prohibitive. 

A practical approach to scalable model selection builds on two highly successful machine learning techniques:  libraries of candidate model nonlinearities and sparse regression.  In particular, the library allows us to choose from a vast range of functional interactions in the network structure.  
Indeed, choosing a large library of functional forms provides the basis to select the nonlinear interactions (models) among nodes of the network.  
By using sparsity promoting techniques like LASSO regression~\cite{Tibshirani:1996} or elastic-net regularization~\cite{Zou2005}, a sparse subset of the components in the library are selected. The information theoretic framework selects nonlinear models that are closest in distance statistically to the empirically measured data. The SINDy method builds, via the function library, an extremely large set of potential models and selects the appropriate terms via sparse regression.  This method has been recently used to successfully infer nonlinear dynamical systems~\cite{Brunton2016pnas}.

\subsection{Contribution of this work}

In this work, we extend the sparse identification of nonlinear dynamical systems (SINDy) algorithm to include rational function nonlinearities in the dynamics. 
Nonlinear dynamics of metabolic and regulatory networks often include rational terms, motivating this innovation.  
It is difficult to construct a library containing rational nonlinearities for use with the sparse regression, since a generic rational function is not simply a sparse linear combination of a small set of basis functions. 
Instead, we write the system as an implicit differential equation in terms of the state and derivatives, and then search for the sparsest vector in the null space of all mixed state and derivative terms.  
For selecting the model terms in implicit-SINDy, we make use of an optimization formulation by Wright et al.~\cite{Wright2009ieeetpami}, and an algorithm using the alternating directions method~\cite{Qu2014}, to find the sparsest vector in the null space.   This selects the active terms in the dynamics.
We demonstrate the algorithm to be robust, accurate and fast when applied to three canonical models of biological networks:  Michaelis-Menten enzyme kinetics, the regulatory network for competence in bacteria, and the metabolic network for yeast glycolysis.

In the following sections, we provide background on the existing SINDy method, and describe the updated algorithm. 
We next validate the algorithm on simulated data for three important biological models: the most fundamental model for enzyme kinetics, a canonical model for regulation of cell differentiation, and a seven-node metabolic network describing glycolysis. 
Finally we discuss the practical application of implicit-SINDy to real biological systems, including overcoming challenges such as noise and increased system size and incorporating perturbative measurements.

\section{Sparse regression for dynamical systems}

Our network inference method is enabled by key innovations around the broadly applicable machine learning methods of sparse regression and overcomplete libraries.  When combined with dynamical systems theory, our robust mathematical architecture is achieved.  

\subsection{LASSO and sparse approximation}
The sparsity promoting $\ell_1$ norm is widely used in a variety of mathematical, scientific, and engineering applications.  Its popularity stems from the fundamental observation that no matter the size of the data from a complex system, there often exists an underlying low-dimensional or {\em sparse} representation of the patterns of interest.  The $\ell_1$ norm is used in fields such as video processing~\cite{Shi:2014}, regression~\cite{Tibshirani:1996}, machine-learning~\cite{Tibshirani:2002,Zou:2006,Wright2009ieeetpami}, signal processing~\cite{Donoho:2006,Candes:2006,Baraniuk:2007}, dynamical systems~\cite{Wang:2011,Ozolicnvs2013pnas,Schaeffer2013pnas,Bright:2013,Brunton2014siads,mackey2014compressive,Brunton2015amr}.

\begin{figure*}[!t]
	\centering
	\includegraphics[width = \textwidth]{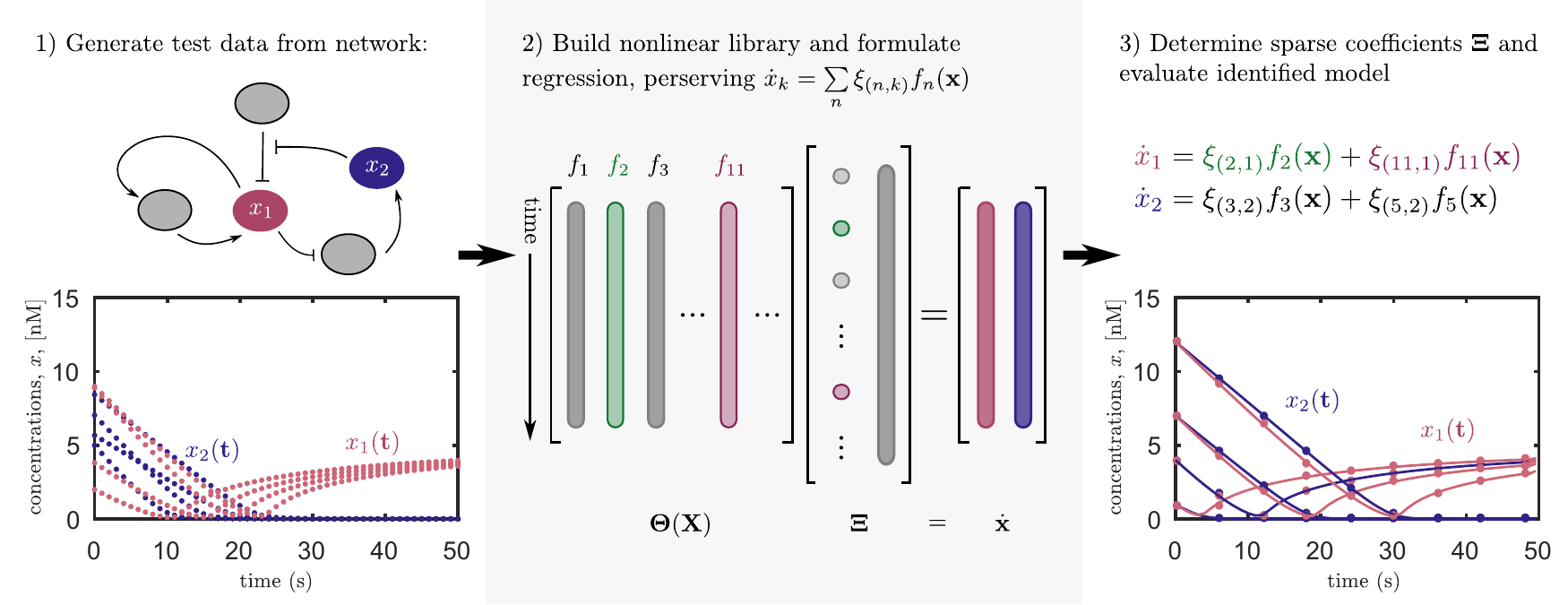}
	\vspace{-.3in}
	\caption{Methodology for sparse identification of nonlinear dynamics (SINDy) from data. First, data is generated from a dynamical system, in this case a biological network. The time series of data is synthesized into a nonlinear function library, and the terms in this library are related to the time derivative by an overdetermined linear regression problem. Enforcing sparsity ensures that only a small number of coefficients are nonzero, identifying the few active terms in the dynamics that are needed to model the system.}
	\label{fig:overview}
	\vspace{-.1in}
\end{figure*}

The least absolute shrinkage and selection operator (LASSO)~\cite{Tibshirani:1996} is one of the earliest sparsity-promoting regression techniques, and it has been successful for feature selection.  
The LASSO solves the linear system of equations 
\begin{eqnarray}
\mathbf{y} = \boldsymbol{\Theta}\mathbf{a},
\end{eqnarray}
by finding a sparse solution vector $\mathbf{a}$ that select features from the over-determined library $\boldsymbol{\Theta}$.  This gives the sparsity criteria for variable selection.  

Given the success of LASSO regression, many other sparsity promoting innovations based on the $\ell_1$ norm have been developed~\cite{Tibshirani:2002,Zou:2006,Line:2011}, including the commonly used elastic-net regularization \cite{Zou2005}.  Closely related to our aims, {\em sparse approximation} exploits the $\ell_1$ norm in order to approximate a set of measurements from a linear combination of vectors from a library~\cite{Donoho:2001,Donoho:2003,Tropp:2003,Gilbert:2003}. The vectors in a library are called atoms.  The atoms in a library can come from pre-computed training data or known bases such as Fourier or wavelets.  The library is deliberately constructed to be overcomplete, thus, creating an overdetermined set of linear equations.  The $\ell_1$ norm overcomes this challenge in order to find a sparse representation of the measurement in the library~\cite{Tropp:2004,Wright2009ieeetpami}.   Like implicit-SINDy, the sparse approximation framework combines two critical ideas:  library building and sparse regression.

\subsection{Sparse identification of nonlinear dynamics (SINDy)}
Discovering dynamical systems from data is an age old pursuit in mathematical physics.  
Historically, this process relied on a combination of high-quality measurements and expert intuition.  
With growing computational power and vast quantities of data, the \emph{automated} discovery of dynamical systems and governing equations is a relatively recent phenomenon.  
Broadly speaking, these techniques may be classified as \emph{system identification}, where methods from statistics and machine learning are used to identify dynamical systems from data.  

Nearly all methods of system identification involve some form of regression of data onto dynamics, and the main distinction between the various techniques is the degree to which this regression is constrained.  
For example, the dynamic mode decomposition (DMD)~\cite{Kutz2016siam} generates best-fit linear models.  
In some cases it is possible to extend these linear regression techniques to nonlinear systems using the Koopman operator~\cite{Koopman1931pnas}, which is an infinite dimensional linear operator acting on functions of observables~\cite{Mezic2013arfm,Budivsic2012chaos,Mezic2005nd}.  
Recent nonlinear regression techniques have produced nonlinear dynamic models that preserve physical constraints, such as conservation of energy~\cite{Majda2012nonlinearity}.  
Genetic programming has also been used to discover dynamical systems and conservation laws from data~\cite{Bongard2007pnas,Schmidt2009science}. 
These methods are highly flexible and impose very few constraints on the form of the dynamics identified.  

Here we review the recent sparse identification of nonlinear dynamics (SINDy) method, which leverages advances in machine learning and sparse regression to discover nonlinear dynamical systems from data~\cite{Brunton2016pnas}.  
SINDy uses sparse regression~\cite{Tibshirani:1996} for improved numerical robustness in noisy overdetermined problems, as opposed to earlier methods~\cite{Wang:2011} that use compressed sensing~\cite{Candes:2006,Donoho:2006,Baraniuk:2007}.

SINDy relies on the fact that many dynamical systems
\begin{eqnarray}
\dot{\mathbf{x}} = \mathbf{f}(\mathbf{x})\label{Eq:NLDyn}
\end{eqnarray}
are sparse in a given function space.  
The relevant terms that are active in the dynamics are solved for using an $\ell_1$-regularized regression that penalizes the number of active terms.  
The general framework for SINDy is shown in Fig.~\ref{fig:overview}.

Algorithmically, time-series data is collected from Eq.~\eqref{Eq:NLDyn}, resulting in a data matrix:
\begin{eqnarray}
\mathbf{X} = \begin{bmatrix} \mathbf{x}(t_1) & \mathbf{x}(t_2) & \cdots & \mathbf{x}(t_m)\end{bmatrix}^T, \label{eq:xmat}
\end{eqnarray} 
where $^T$ denotes the matrix transpose. 
The matrix $\mathbf{X}$ is $m\times n$, where $n$ is the dimension of the state $\mathbf{x}\in\mathbb{R}^n$ and $m$ is the number of measurements of the state in time.  
For our purposes the state variables are the measured biological components in the network (enzymes, metabolites, transcription factors etc.).  
Similarly, the matrix of derivatives
\begin{eqnarray}
\dot{\mathbf{X}} = \begin{bmatrix} \dot{\mathbf{x}}(t_1) &  \dot{\mathbf{x}}(t_2) & \cdots &  \dot{\mathbf{x}}(t_m)\end{bmatrix}^T,
\end{eqnarray} 
is collected or computed from the state data in $\mathbf{X}$; the total-variation regularized~\cite{Rudin1992physd} derivative~\cite{Chartrand2011isrnam} provides numerically robust method to compute derivatives from noisy data.  

Next, a library of candidate nonlinear functions is constructed from $\mathbf{X}$:
\begin{eqnarray}
\boldsymbol{\Theta}(\mathbf{X}) = \begin{bmatrix} \mathbf{1} & \mathbf{X} & \mathbf{X}^2 & \cdots & \mathbf{X}^d  & \cdots &   \sin(\mathbf{X}) & \cdots  \end{bmatrix},\label{Eq:NLLibrary}
\end{eqnarray}
where $\mathbf{X}^d$ denotes the matrix containing all possible column vectors obtained from time-series of the $d$-th degree polynomials in the state vector $\mathbf{x}$.  For example, for a system with two states $\mathbf{x} = \begin{bmatrix} x_1, & x_2\end{bmatrix}^T$, the matrix $\mathbf{X}^2 = \begin{bmatrix} x_1^2(\mathbf{t}), & (x_1x_2)(\mathbf{t}), & x_2^2(\mathbf{t})\end{bmatrix}$, where $\mathbf{t}$ is a vector of times at which the state is measured.  

It is now possible to relate the time derivatives in $\dot{\mathbf{X}}$ to the candidate nonlinearities in $\boldsymbol{\Theta}(\mathbf{X})$ by:
\begin{eqnarray}
\dot{\mathbf{X}} = \boldsymbol{\Theta}(\mathbf{X})\boldsymbol{\Xi},\label{Eq:SINDy1}
\end{eqnarray}
where each column $\boldsymbol{\xi}_k$ in $\boldsymbol{\Xi}$ is a vector of coefficients that determines which terms are active in the $k$-th row equation of Eq.~\eqref{Eq:NLDyn}.  
To enforce sparsity in the dynamics, we solve for each column $\boldsymbol{\xi}_k$ using sparse regression, such as the LASSO~\cite{Tibshirani:1996}:
\begin{eqnarray}
\boldsymbol{\xi}_k = \text{argmin}_{\boldsymbol{\xi}_k'}\|\dot{\mathbf{X}}_k - \boldsymbol{\Theta}(\mathbf{X})\boldsymbol{\xi}_k'\|_2+\lambda \|\boldsymbol{\xi}_k'\|_1,
\end{eqnarray}
where $\dot{\mathbf{X}}_k$ is the $k$-th column of $\dot{\mathbf{X}}$.  Once the sparse coefficient vectors $\boldsymbol{\xi}_k$ are determined, a model of the nonlinear dynamical system may be constructed:
\begin{eqnarray}
\dot{\mathbf{x}} = \boldsymbol{\Xi}^T\left(\boldsymbol{\Theta}(\mathbf{x}^T)\right)^T.
\end{eqnarray}

Using sparse regression to identify active terms in the dynamics from the candidate library $\boldsymbol{\Theta}(\mathbf{X})$ is a convex optimization.  
The alternative is to apply a separate constrained regression on every possible subset of nonlinearities, and then to choose the model that is both accurate and sparse.   
This brute-force search is intractable, and the SINDy method makes it possible to select the sparse model in this combinatorially large set of candidate models. 

\begin{figure*}[!t]
	\centering
	\includegraphics[width=\textwidth]{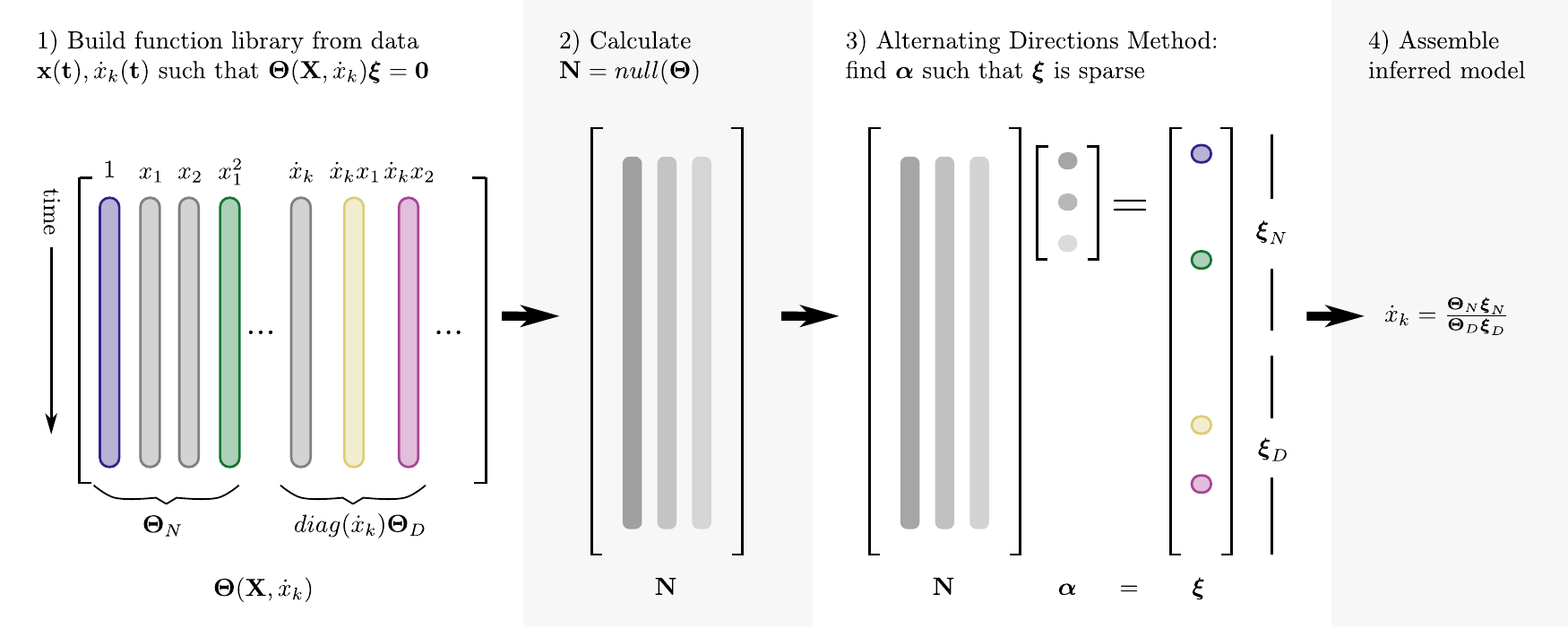}
	\vspace{-.3in}
	\caption{SINDy algorithm for rational functions. Assemble a matrix, $\boldsymbol{\Theta}(\mathbf{X},\dot{x}_k)$ where each column is a nonlinear function evaluated for time series data $x_1, x_2, x_3 .. $ and one species derivative $\dot{x}_k(\mathbf{t})$. Next, calculate  $\mathbf{N}$ a orthonormal basis for the null space of $\boldsymbol{\Theta}$. Then, use an alternating directions method \cite{Qu2014} to find a sparse vector, $\boldsymbol{\xi}$, in the null space. The sparse vector $\boldsymbol{\xi}$, then satisfies $\boldsymbol{\Theta} \boldsymbol{\xi} = 0$. Using the sparse coefficients from $\boldsymbol{\xi}$ and the functional library $\boldsymbol{\Theta}$ assemble the inferred model. This algorithm must be performed for the derivative $\dot{x}_k$ of each species. }
	\label{fig:SINDy2}
	\vspace{-.1in}
\end{figure*}

The polynomial and trigonometric nonlinearities in Eq.~\eqref{Eq:NLLibrary} are sufficient for a large class of dynamical systems.  
For example, evaluating all polynomials up to order $n$ is equivalent to assuming that the biological network has dynamics determined by mass action kinetics up to $n$-mers (monomers, dimers, trimers, etc.). 
However, if there are time-scale separations in the mass action kinetics, fast reactions are effectively at steady state, and the remaining equations contain rational functions \cite{Gunawardena2014}.
As we consider systems where the dynamics include rational functions, constructing a comprehensive library becomes more complicated.
If we generate all rational nonlinearities:
\begin{eqnarray}
f(\mathbf{x}) = \frac{f_N(\mathbf{x})}{f_D(\mathbf{x})}, \label{Eq:ratfunc}
\end{eqnarray}
where $f_N(\mathbf{x})$ and $f_D(\mathbf{x})$ are both polynomial functions, the library would be prohibitively large for computational purposes.
Therefore, we develop a computationally tractable framework in the next section for functional library construction that accounts for dynamics with rational functions.

\section{Inferring nonlinear dynamical systems with rational functions}
Many relevant dynamical systems contain rational functions in the dynamics, motivating the need to generalize the SINDy algorithm to include more general nonlinearities than simple polynomial or trigonometric functions.  
The original SINDy algorithm bypasses the computation and evaluation of all $N_p$ candidate regression models, as enumerated in Eq.~\eqref{eq:nummodels}, by performing a sparse approximation of the dynamics in a library constructed from the $N_m$ candidate monomial features. 
However, it is not possible to simply apply the original SINDy procedure and include rational functions, since generic rational nonlinearities are not sparse linear combinations of a small number of rational functions.  
Instead, it is necessary to modify the sparse dynamic regression problem to solve for the sparsest \emph{implicit} ordinary differential equation according to the following procedure.

Consider a dynamical system of the form in Eq.~\eqref{Eq:NLDyn}, but where the dynamics of each $k=1,2,\cdots, n$ variables may contain rational functions:
\begin{eqnarray}
\dot{x}_k = \frac{f_N(\mathbf{x})}{f_D(\mathbf{x})}
\end{eqnarray}
where $f_N(\mathbf{x})$ and $f_D(\mathbf{x})$ represent numerator and denominator polynomials in the state variable $\mathbf{x}$. 
For each equation, it is possible to multiply both sides by the denominator polynomial, resulting in the equation:
\begin{eqnarray}
f_N(\mathbf{x})-f_D(\mathbf{x}) \dot{x}_k  = 0. \label{eq:NDzero}
\end{eqnarray}

The implicit form of Eq.~\eqref{eq:NDzero} motivates a generalization of the function library $\boldsymbol{\Theta}$ in Eq.~\eqref{Eq:NLLibrary} in terms of the state $\mathbf{x}$ and the derivative $\dot{x}_k$:
\begin{eqnarray}
\boldsymbol{\Theta}(\mathbf{X},\dot{x}_k(\mathbf{t})) = \begin{bmatrix} \boldsymbol{\Theta}_N(\mathbf{X}) & \text{diag}\left(\dot{x}_k(\mathbf{t})\right)\boldsymbol{\Theta}_D(\mathbf{X})\end{bmatrix}.\label{Eq:AugmentedLibrary}
\end{eqnarray}
The first term, $\boldsymbol{\Theta}_N(\mathbf{X})$, is the library of numerator monomials in $\mathbf{x}$, as in Eq.~\eqref{Eq:NLLibrary}.  
The second term, $\text{diag}\left(\dot{x}_k(\mathbf{t})\right)\boldsymbol{\Theta}_D(\mathbf{X})$, is obtained by multiplying each column of the library of denominator polynomials $\boldsymbol{\Theta}_D(\mathbf{X})$ with the vector $\dot{x}_k(\mathbf{t})$ in an element-wise fashion. For a single variable $x_k$, this would give the following:
\begin{eqnarray}
\text{diag}(\dot{x}_k(\mathbf{t}))\boldsymbol{\Theta}(\mathbf{X}) \!=\! \begin{bmatrix} \dot{x}_k(\mathbf{t}) \,\, (\dot{x}_k x_k)(\mathbf{t}) \,\, (\dot{x}_k x_k^2)(\mathbf{t}) \,\, \dotsc \,\, \end{bmatrix}.
\end{eqnarray}

A schematic of this library is shown in Fig.~\ref{fig:SINDy2}. 
In most cases, we will use the same polynomial degree for both the numerator and denominator library, so that $\boldsymbol{\Theta}_N(\mathbf{X})=\boldsymbol{\Theta}_D(\mathbf{X})$.  
Thus, the augmented library in Eq.~\eqref{Eq:AugmentedLibrary} is only twice the size of the original polynomial library in Eq.~\eqref{Eq:NLLibrary}.  

\begin{figure}[t!]
	\centering
	\vspace{-.05in}
	\includegraphics[width=0.5\textwidth]{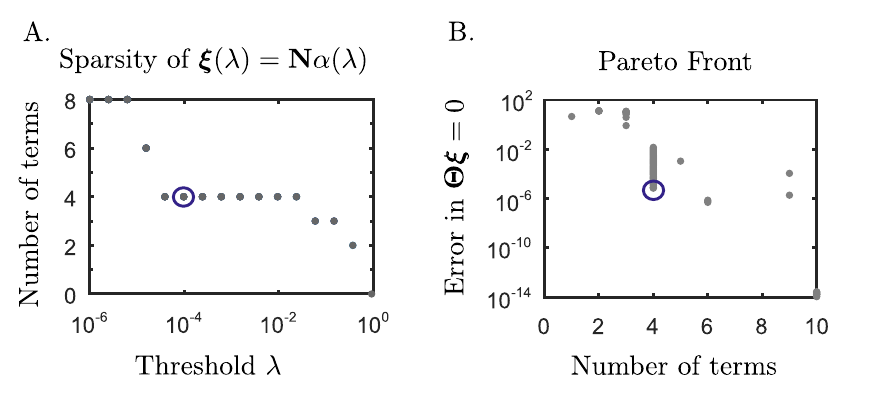}
	\vspace{-.3in}
	\caption{A. Increasing the sparsity threshold $\lambda$ during ADM creates coefficient vectors, $\boldsymbol{\xi}$, with monotonically decreasing number of terms. B. For each $\boldsymbol{\xi}(\lambda)$ we calculate an error as $|\boldsymbol{\Theta}\boldsymbol{\xi}|$, and produce the Pareto Front. For the cases tested here, a large cliff in the error indicates the best choice of $\boldsymbol{\xi}(\lambda)$ (circled on A. and B.) and the most parsimonious model.  }
	\label{fig:sparsity}
	\vspace{-.15in}
\end{figure}

We may now write the dynamics in Eq.~\eqref{eq:NDzero} in terms of the augmented library in Eq.~\eqref{Eq:AugmentedLibrary}:
\begin{equation}
\boldsymbol{\Theta}(\mathbf{X},\dot{x}_k(\mathbf{t}))\boldsymbol{\xi}_k = \mathbf{0}. \label{Eq:Sindy2}
\end{equation}
The sparse vector of coefficients $\boldsymbol{\xi}_k$ will have non-zero entries for the terms active in the nonlinear dynamics.  
However, it is not possible to use the same method of sparse regression as in the original SINDy algorithm, i.e. to find the sparsest vector $\boldsymbol{\xi}_k$ that satisfies Eq.~\eqref{Eq:Sindy2}, since the sparsest vector would be identically zero.  

To find the sparsest non-zero vector $\boldsymbol{\xi}_k$ that satisfies Eq.~\eqref{Eq:Sindy2}, we note that any such vector will be in the null space of $\boldsymbol{\Theta}$.  
After identifying the null space of $\boldsymbol{\Theta}$, we need only find the sparsest vector in this subspace.  
Although this is a non-convex problem, there are straightforward algorithms based on the alternating directions method (ADM) developed by Qu \emph{et al.}~\cite{Qu2014} to identify the sparsest vector in a subspace.  

\subsection{Algorithm for sparse selection of rational functions.}

The algorithm for finding $\boldsymbol{\xi}_k$ is as follows.  First, we build our functional library $\boldsymbol{\Theta}(\mathbf{X},\dot{x}_k(\mathbf{t}))$ using both the time series data of the state variables and derivative, as discussed above.
Second, we calculate a matrix, $\mathbf{N}$, with columns spanning the null space of $\boldsymbol{\Theta}$. 
We wish to find the linear combinations of columns in $\mathbf{N}$ that produces a sparse vector $\boldsymbol{\xi}$.
For this third step, we use the alternating directions method developed by Qu \emph{et al.} \cite{Qu2014} that finds the sparsest vector in a subspace.
We enforce some magnitude of sparsity using a threshold, $\lambda$.
For the fourth and final step, we select the active nonlinear functions using $\boldsymbol{\xi}$ and $\boldsymbol{\Theta}$, and assemble the inferred model.

As the appropriate $\lambda$ is unknown {\sl a priori}, we repeat the third and fourth steps for varying $\lambda$.
Increasing $\lambda$ increases the sparsity (decreasing the number of terms) in $\boldsymbol{\xi}$, as shown in Fig.~\ref{fig:sparsity}A.
Each $\boldsymbol{\xi}(\lambda)$ produces an inferred model of varying accuracy and sparsity. 
From these models we calculate a Pareto front and select the most parsimonious model, as shown in Fig.~\ref{fig:sparsity}B.
The most parsimonious model is readily identifiable at the sharp drop-off on the Pareto plot.
As we will show, this method succeeds at identifying the correct rational terms and coefficients.

\subsection{General formulation for implicit ODEs}
The procedure above may be applied to identify more general implicit ordinary differential equations, beyond those just containing rational function nonlinearities.  
The library $\boldsymbol{\Theta}(\mathbf{X},\dot{x}_k(\mathbf{t}))$ contains a subset of the columns of the library $\boldsymbol{\Theta}(\begin{bmatrix} \mathbf{X} & \dot{\mathbf{X}}\end{bmatrix})$, which is obtained by building nonlinear functions of the state $\mathbf{x}$ and derivative $\dot{\mathbf{x}}$.  

\begin{figure*}[!t]
	\centering
	\includegraphics[width=\textwidth]{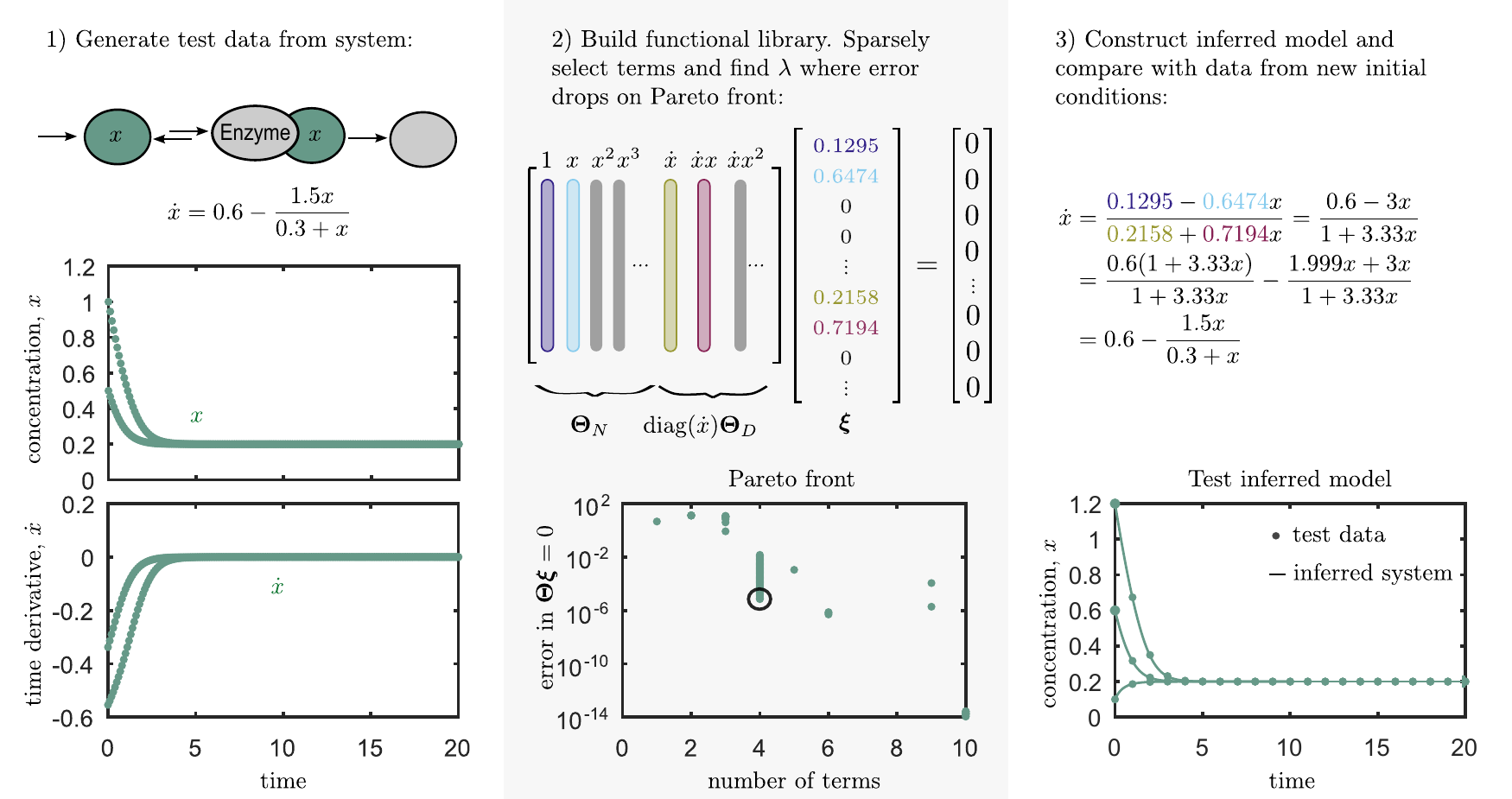}
	\vspace{-.3in}
	\caption{Algorithm applied to the Michaelis-Menten kinetics for an enzymatic reaction. Step 1) Generate two time series of the single state variable, $x(\mathbf{t})$, and time derivative, $\dot{x}(\mathbf{t})$. Step 2) Discovered active functions and their corresponding coefficients are indicated by color. The error drops sharply at 4 terms on the Pareto front (circled) . The most parsimonious model has four active functions: two in the numerator and two in the denominator (indicated by color). Step 3) Allowing for rational function factorization, the inferred model is equivalent to the original model. }
	\label{fig:MM}
\end{figure*}

Identifying the sparsest vector in the null space of $\boldsymbol{\Theta}(\begin{bmatrix} \mathbf{X} & \dot{\mathbf{X}}\end{bmatrix})$ provides more flexibility in identifying nonlinear equations with mixed terms containing various powers of any combination of derivatives and states. 
For example, the system given by
\begin{eqnarray}
\dot{x}^3x - \dot{x}x^2 - x^3 = 0
\end{eqnarray}
may be encoded as a sparse vector in the null space of $\boldsymbol{\Theta}(\begin{bmatrix} \mathbf{X} & \dot{\mathbf{X}}\end{bmatrix})$.  
It is also straightforward to extend the formulation to include higher order derivatives, by increasing the features in the $\boldsymbol{\Theta}$ library.  
For example, second-order implicit dynamical systems may be formulated in the following library:
\begin{eqnarray}
\boldsymbol{\Theta}\left(\begin{bmatrix}\mathbf{X} & \dot{\mathbf{X}} & \ddot{\mathbf{X}}\end{bmatrix}\right).
\end{eqnarray}
The generality of this approach enables the identification of many more systems of interest, in addition to those systems with rational function nonlinearities explored below.  

\section{Results}

The implicit-SINDy architecture is tested on a number of canonical models of biological networked dynamical systems.  Validation of the method on these models allows for potential broader application.  We demonstrate that the method is fast, accurate and robust for inferring Michaelis-Menten enzyme kinetics, the regulatory network for competence in bacteria, and the metabolic network for yeast glycolysis.

\subsection{Simple example: Michaelis-Menten kinetics}

Perhaps the most well known model for enzyme kinetics is the Michaelis-Menten model \cite{MichaelisL&MentenM1913, Johnson2011}. This model captures the dynamics of an enzyme binding and unbinding with a substrate ($x$), and then reacting irreversibly to produce a product, as shown in Fig.~\ref{fig:MM}. A separation of time-scales argument, where binding and unbinding dynamics are fast, or a more general steady state assumption~\cite{Briggs1925}, reduces the dynamics to a single state-variable equation with a rational function in the dynamics. Traditionally, biochemists vary the initial concentration of $x$ in a titration experiment to fit the Michaelis-Menten equation to the data. 

Using time series data from only two initial concentrations, our algorithm extracts the correct functional form from a larger subset of possible functions and fits the coefficients accurately (Fig.~\ref{fig:MM}). First we generate data from the single dynamic equation
\begin{equation}
	\dot{x} = j_x - \frac{V_{max} x}{K_m + x},
\end{equation}
with some flux source of $x$, $j_x$, and an enzymatic reaction of the Michaelis-Menten form consuming $x$. Here, $V_{max}$ is the maximum rate of the reaction and $K_m$ is the concentration of half-maximal reaction rate. Generally the time series data of the concentration, $x(\mathbf{t})$, is measurable, while the time series data for the derivative can be calculated from $x(\mathbf{t})$. 

Next, we apply implicit-SINDy to determine the coefficient vector $\boldsymbol{\xi}$ and sparsely select the active functions in the dynamics. The library contains polynomial terms up to degree four and has 10 columns. 
The Pareto front for this system has a sharp drop off in error from around $0.01$ to $10^{-5}$ at four terms, indicating the $\lambda$ for the most parsimonious model. The associated $\boldsymbol{\xi}$ selects 4 active terms from the function library. 
\begin{figure*}[!t]
	\centering
	\includegraphics[width=\textwidth]{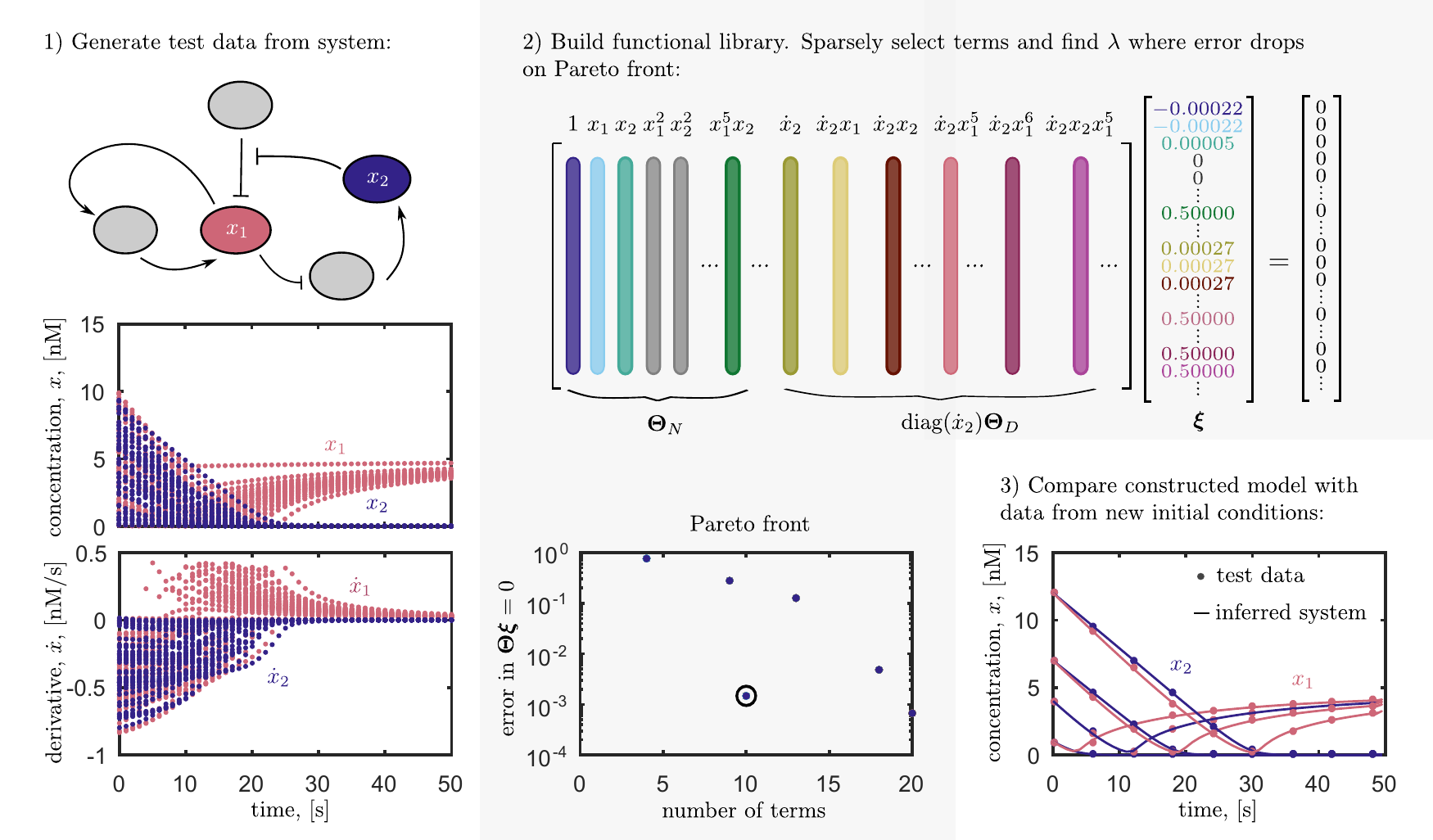}
	\vspace{-.15in}
	\caption{Algorithm applied to a regulatory network with two measured state-variables: $x_1$ (blue) and $x_2$ (pink). Pointed arrows indicate activation and blunted arrows indicate repression. The functional library $\boldsymbol{\Theta}$, sparse vector $\boldsymbol{\xi}$, and Pareto Front for $x_2$ are shown in Step 2. 10 of 56 terms are active in the library for the most parsimonious model: 4 in the numerator and 6 in the denominator. The inferred model is nearly equivalent to the original model for test data. }
	\label{fig:reg}
\end{figure*}
Finally, the rational function constructed from $\boldsymbol{\Theta}$ and $\boldsymbol{\xi}$ needs to be factored to be interpreted as the source flux and Michaelis-Menten terms. When rearranged, the coefficients match the original system. Unsurprisingly, the inferred model matches the original model for time series generated from new initial conditions that were not used in the training data.

\subsection{Regulatory network: {\sl B. subtilis} competence}

Having shown that our method works for the simplest rational model relevant to biological networks, we next test it on a regulatory model with two state variables \cite{Suel2006}. S\"{u}el et al.~\cite{Suel2006} demonstrated that a dynamic gene network enables cells to switch between multiple behaviors -- in this case {\sl B. subtilis} bacteria switch between taking up DNA from the environment (competence) and vegetative growth. 
Other regulatory networks such as the circadian clock \cite{Hardin2005, Leloup2003} and cell cycle oscillators have similar structure and dynamics. In particular, similar dynamics may drive cancer-relevant systems like the tumor suppressor p53 \cite{Batchelor2009}. 

The dynamics of regulatory system with two states can be described by the following two non-dimensional equations:
\begin{subequations}
	\begin{eqnarray}
		\dot{x}_1 &=& a_1 + \frac{a_2 x_1^2}{a_3 + x_1^2} - \frac{x_1}{1+ x_1 + x_2},\label{eq:x1reg}\\ 
		\dot{x}_2 &=& \frac{b_1}{1+ b_2 x_1^5} - \frac{x_2}{1 + x_1+x_2 }. \label{eq:x2reg}
	\end{eqnarray}
\end{subequations}
These two equations are a reduction of dynamical system with six states. Each rational function arises from a steady state (or time-scale separation) assumption about the regulatory processes. The second term (scaled by $a_2$) in Eq.~\eqref{eq:x1reg} represents protein $x_1$, ComK, activating its own production in an autoregulatory, positive feedback loop. The first term (scaled by $b_1$) in Eq.~\eqref{eq:x2reg}, describes $x_1$ mediated repression of $x_2$, ComS, in a negative feedback loop. Both of these terms have a Hill-function form, where the power indicates the number of $x_1$ proteins involved cooperatively in the regulatory complex \cite{Alon2007a}. The combination of positive and negative feedback results in the network's functional capabilities. The last term in Eqs.~\eqref{eq:x1reg} and~\eqref{eq:x2reg} describes degradation of $x_1$ and $x_2$, mediated by a third unmeasured protein, MecA.

Using this model we generate 40 time series for the regulatory system, as shown in Fig.~\ref{fig:reg}.  This model challenges our method in two ways. First, the method must correctly identify the dynamic dependence on two state variables. Second, the model contains polynomial functions up to the 5th degree in the numerator of one term in Eq.~\eqref{eq:x2reg}. To include this term, the library must contain polynomials up to degree six. Even without knowing the highest polynomial power ahead of time, it is possible to use the implicit-SINDy by trying libraries of increasing polynomial degree. If the library does not have all the required terms, there will be no clear drop off in the Pareto front as there is in Fig.~\ref{fig:reg}. 
\begin{figure*}[!t]
	\centering
	\includegraphics[width = \textwidth]{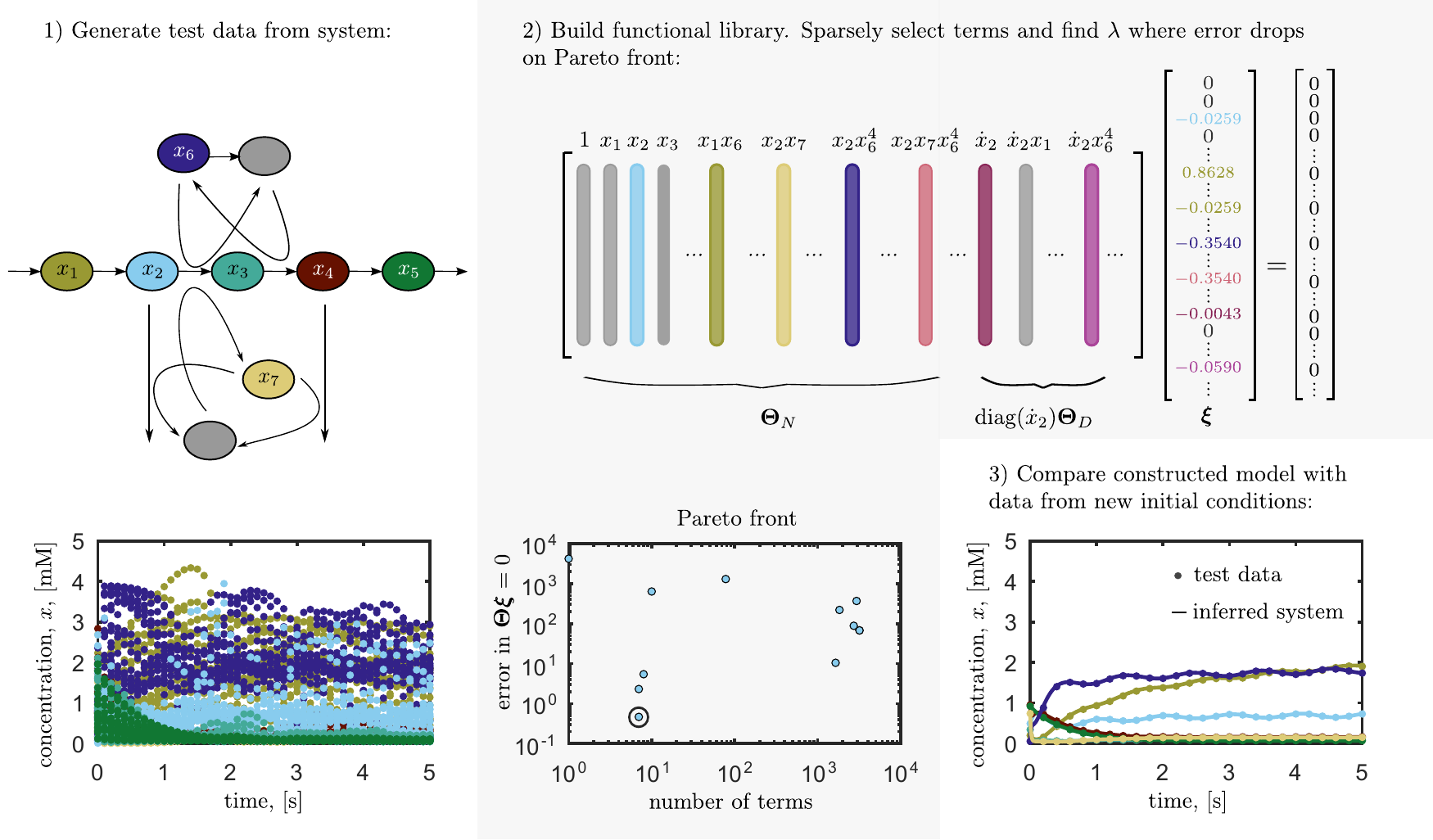}
	\vspace{-.15in}
	\caption{Algorithm applied to a metabolic network for yeast glycolysis. Step 1). The 7 measured state variables are indicated with separate colors consistent with time series. Arrows indicate reactions between components. Not all data required to infer the network is plotted. Step 2) shows functional library $\boldsymbol{\Theta}$ and corresponding coefficient vector $\boldsymbol{\xi}$ for $x_2$.  Seven functions are active: 5 in the numerator and 2 in the denominator. Step 3) The inferred model is nearly equivalent to the original system for test data. One time series is shown for all 7 state variables. }
	\label{fig:met}
\end{figure*}
The library with degree six polynomials in two state variables contains 56 columns, of which 10 are active in the most parsimonious model for $x_2$ dynamics. The constructed models for $x_1$ and $x_2$ match almost exactly with test data generated from the original model. As with our first example, the extracted rational function can be factored to recover exactly the form of Eq.~\eqref{eq:x2reg}.  Additionally, the coefficients identified are within 2\% error of the true coefficients shown in Table~\ref{tab:reg}. 

\begin{table}[!h]
	\caption{Parameter identification for regulatory network} \label{tab:reg}
	\center
	\vspace{-3mm}
	\begin{tabular}{|l|l|l|l|}
		\hline
		Parameter  & units  & True Value   	   & Extracted value \\ \hline
		$a_1$      & [nM/s] & 0.004         & 0.00393              \\
		$a_2$  		& [nM/s] & 0.07        & 0.07006             \\
		$a_3$  & [nM]       & 0.04       	   & 0.04000             \\ \hline
		$b_1$  & [nM/s]     & 0.82          & 0.8148             \\
		$b_2$  & [nM]       & 1854.5 	  	   & 1851.9             \\ \hline
	\end{tabular}
\end{table}

\subsection{Metabolic network: yeast glycolysis }

As a final example, we test our method on a metabolic network. Glycolysis, the process of breaking down glucose to extract energy (ATP and NADPH), is part of central metabolism for all cells. Uncovering the metabolic network for glycolysis took over 100 years from its initial discovery by Pasteur \cite{Racker1980}. Accelerated inference of metabolic networks would aid metabolic disease intervention \cite{Yarmush2003}. Bacteria perform a wide range of yet-to-be discovered chemistry which could be harnessed through metabolic engineering to produce high-value products such as drugs and biofuels \cite{Ellis2012a}\cite{Ducat2011}.

Not only does the yeast glycolysis model we analyze have a larger number of interacting state variables, but it is  oscillatory\cite{Wolf2000}. This network has also been previously analyzed as a test case for model inference \cite{Schmidt2011pb}. The network shown in Fig.~\ref{fig:met} has three equations with rational functions and four with polynomials:
\begin{subequations}
	\begin{eqnarray}
		\dot{x}_1 &=& c_1 + \frac{c_2 x_1 x_6}{1+ c_3 x_6^4}\\ \label{eq:x1met}
		\dot{x}_2 &=& \frac{d_1 x_1 x_6}{1+ d_2 x_6^4} + d_3 x_2  - d_4 x_2 x_7 \\ \label{eq:x2met}
		\dot{x}_3 &=& e_1 x_2 + e_2 x_3+ e_3 x_2 x_7 + e_4 x_3 x_6  \\ \label{eq:x3met}
		\dot{x}_4 &=& f_1 x_3 + f_2 x_4 + f_3 x_5+ f_4 x_3 x_6 +f_5 x_4 x_7  \\
		\dot{x}_5 &=& g_1 x_4 + g_2 x_5 \\
		\dot{x}_6 &=& \frac{h_1 x_1 x_6}{1+ h_2 x_6^4} + h_3 x_3 +h_5 x_6  + h_4 x_3 x_7\\ \label{eq:x6}
		\dot{x}_7 &=& j_1 x_2 +j_2 x_2 x_7 +j_3 x_4 x_7.
	\end{eqnarray}
\end{subequations}
Given sufficient data, implicit-SINDy correctly infers the network structure and coefficients. In Fig.~\ref{fig:met}, we show the sparsely selected terms and Pareto front for Eq.~\eqref{eq:x2met}. The method correctly selects 7 terms from a library of 3432 functions: 5 in the numerator and 2 in the denominator. Table~\ref{tab:met} shows the true and extracted coefficient values for the model. Some of the parameters, $c_1$ for example, had incorrect functional dependence on $x_6$ after factoring the discovered polynomial. However these dependencies were very small ($<0.1\%$ error). 

\begin{table}[!t]
	\caption{Parameter identification for metabolic network} \label{tab:met}
	\center
	\vspace{-3mm}
	\begin{tabular}{|l|l|l|l|}
		\hline
		Parameter	& units    			& True Value 	& Extracted value \\ \hline
		$c_1$   	& [mM/min]  		& 2.5         	& $2.5002^*$         \\
		$c_2$       & [1/(mM min)]		& -100        	& -99.7979        \\
		$c_3$       & [mM${}^{-1/4}$] 	& 13.6769     	& 13.6489         \\ \hline
		$d_1$       & [1/(mM min)]	 	& 200         	& 200.6512        \\
		$d_2$       & [mM${}^{-1/4}$]	&13.6769      	& 13.7209               \\
		$d_3$       & [1/min]			& -6          	& $-5.9998^*$     \\
		$d_4$       & [1/(mM min)]		& -6          	& $-5.9998^*$    \\ \hline
		$e_1$       & [1/min]			& 6           	& 6.0133         \\
		$e_2$       & [1/min]			& -64         	& -64.140         \\
		$e_3$       & [1/(mM min)]		& 6           	& -6.0133          \\
		$e_4$       & [1/(mM min)]		& 16          	& 16.0333         \\ \hline
		$f_1$       & [1/min]			& 64          	&  63.6145        \\
		$f_2$       & [1/min]			& -13         	& -12.9277        \\
		$f_3$       & [1/min]			& 13          	&  12.9277        \\
		$f_4$       & [1/(mM min)]		& -16         	& -15.9036               \\
		$f_5$       & [1/(mM min)] 		& -100        	& -99.3976               \\ \hline
		$g_1$       & [1/min] 			& 1.3         	& 1.3002               \\
		$g_2$       & [1/min]	 		& -3.1        	& -3.1003               \\ \hline
		$h_1$       & [1/(mM min)]		& -200        	& $-200.0^\dagger$     \\
		$h_2$       & [mM${}^{-1/4}$]   & 13.6769     	& $13.6769^\dagger$    \\
		$h_3$       & [1/min] 			& 128         	& $128.0^\dagger$      \\
		$h_4$       & [1/(mM min)]      & -1.28       	& $-1.280^\dagger$     \\
		$h_5$       & [1/min] 		    & -32       	& $-32.00^\dagger$     \\ \hline
		$j_1$       & [1/min]			& 6           	& 6.0102               \\
		$j_2$       & [1/(mM min)] 		& -18         	& -18.0408               \\
		$j_3$       & [1/(mM min)] 		& -100        	& -100.2449               \\ \hline
	\end{tabular}
	\\
	\vspace{+2mm}
	\begin{flushleft}
		Values with a ${}^*$ have errors with functional dependence on $x_6$. We show the leading order term assuming a Taylor series expansion. ${}^\dagger$ Eq.~\eqref{eq:x6} required more data to identify, and this allowed extraction of the coefficients to a higher precision.
	\end{flushleft}	
\end{table}
The equations for $x_3, \; x_4, \; x_5$, and $x_7$ did not require a library with polynomials up to degree six, and could be inferred more quickly. On the other hand, Eq.~\eqref{eq:x6}, required over twice as many measurements as the other equations with rational functions (Eqs.~\eqref{eq:x2met} and~\eqref{eq:x1met}).

\section{Conclusions}
In this work we developed an implicitly formulated method for sparse identification of nonlinear dynamics: implicit-SINDy.   The method allows for constructing nonlinear dynamics with rational functions using a library of functional forms that is still computationally manageable for reasonably-sized biological networks.  An alternating directions method for selection of a sparse vector in the null space of the library \cite{Qu2014} enables us to find and construct a parsimonious model from the full library.  Using implicit-SINDy on data generated from three biological models (enzyme kinetics, regulation, and metabolism), we are able to accurately reconstruct the underlying system in each case.  Indeed, we correctly recover the coefficients to within 2\% of the original values. These results make implicit-SINDy a promising method for model discovery of biological networks.


%
SINDy is a {\sl data-driven} methodology, meaning it selects the connectivity and dynamics based on the information content of the data alone.
It has many advantages, including the fact that there is no parameter tuning in the inferred models aside from a sparsity threshold which is determined by a Pareto front.
Moreover, implicit-SINDy greatly expands our ability to rapidly select a model from a large class of candidate dynamical systems, even when nonlinear derivative terms are present.  
In practice, the method functions much like a highly efficient unsupervised learning algorithm, sparsely selecting dynamics from a large library of possibilities.
It differs from information theoretic techniques where a number of viable models are posited and selection of the best model is based upon the minimization of information loss.  
Such alternative techniques rely on physical insight (supervised learning) to generate individual models, thus potentially limiting the dynamical systems considered.
Given the large number of biological models driven by mass-action kinetics, the implicit-SINDy method can be a critically enabling method for data-driven discovery of underlying biological principals.


There are many intriguing future directions for the method, both in theory and practice. 
In theory, the connection between the implict-SINDy selection process and information criteria such as AIC and BIC remains an open question.  We hope to pursue this further in order to establish a rigorous statistical connection between information theory metrics and sparse selection.  In practice, two remaining challenges exist for practical implementation:  (i) improving robustness to noise and (ii) reducing the number of time-series measurements.   As one would expect, noise compromises the calculation of the null space of a matrix.  Recent work by Gavish and Donoho provided a general method for recovering a low-rank matrix (and therefore null space) from noisy data~\cite{Gavish2014}.  Such a thresholding technique may be used to make implicit-SINDy more robust to noise.  It remains an open question how long a time-series must be sampled in order to accurately produce the underlying dynamics.  This will be also considered in future work.

\section*{Acknowledgment}

SLB acknowledges support from the U.S. Air Force Center of Excellence on Nature Inspired Flight Technologies and Ideas (FA9550-14-1-0398). 
JLP and NMM would like to thank Bill and Melinda Gates for their active support of the Institute for Disease Modeling and their sponsorship through the Global Good Fund. 

\bibliographystyle{IEEEtran}
\bibliography{biball,extrabib}

\ifCLASSOPTIONcaptionsoff
\newpage
\fi


%
%
%

\end{document}